\documentclass[11pt, psamsfonts]{amsart}
\usepackage{amsmath}
\usepackage{amsthm}
\usepackage{amssymb}
\usepackage{amscd}
\usepackage{amsfonts}
\usepackage{amsbsy}
\usepackage{epsfig}
\usepackage{graphicx}
\usepackage{psfrag}

\newtheorem{theorem}{Theorem}
\newtheorem{remark}{Remark}
\newtheorem{proposition}{Proposition}
\newtheorem{lemma}{Lemma}
\newtheorem{corollary}{Corollary}
\newtheorem*{claim}{Claim}

\newfont\bbf{msbm10 at 12pt}
\def\eps{\varepsilon}
\def\phi{\varphi}
\def\R{{\mathbb R}}
\def\C{{\mathbb C}}
\def\N{{\mathbb N}}

\def\disk{{\mathbb D}}

\def\J{{\mathcal J}}

\def\F{{\mathcal F}}

\def\W{{\mathcal W}}

\def\es{{\emptyset}}
\def\sm{\setminus}

\def\dist{\mbox{dist}}

\def\modulus{\mbox{mod}}

\def\bd{\partial }
\def\le{\leqslant}
\def\ge{\geqslant}

\newcommand{\st}{such that }

\newcommand{\poin}{Poincar\'{e} }

\newcommand{\bcls}{Borel--Cantelli Lemma}

\newcommand{\z}{\zeta}

\begin{document}
\title{Complex maps without invariant densities}
\author{Henk Bruin, Mike Todd}
\thanks{This research was supported by EPSRC grant GR/S91147/01}
\subjclass[2000]{ Primary: 37F35; Secondary 37F10, 37C40, 37D35,
37F25}

\begin{abstract}
\noindent We consider complex polynomials $f(z) = z^\ell+c_1$ for
$\ell \in 2\N$ and $c_1 \in \R$, and find some combinatorial types
and values of $\ell$ such that there is no invariant probability
measure equivalent to conformal measure on the Julia set. This
holds for particular Fibonacci-like and Feigenbaum combinatorial
types when $\ell$ sufficiently large and also for a class of
`long--branched' maps of any critical order.
\end{abstract}

\maketitle

\section{Introduction}

The question of finding invariant probability measures that are
absolutely continuous with respect to a natural reference measure
is central in many areas of dynamical systems. Such invariant
probabilities are called {\em acip}s. In dynamical systems on
manifolds, the reference measure is usually taken to be Lebesgue
(or some Riemannian volume). In this paper we consider unicritical
polynomials $f(z) = z^\ell + c$ on the complex plane with critical
order $\ell \in 2\N$ and parameter $c \in \R$, and in this
setting, the natural reference measure is $\delta$-conformal
measure $\mu_\delta$, supported on the Julia set $\J = \J(f)$.
This probability measure is defined by the relation
\[
\mu_\delta(A) = \int_A |Df(z)|^\delta d\mu_\delta(z),
\]
whenever $f:A \to f(A)$ is one-to-one. There is a
$\delta$--conformal measure for at least one $\delta \in [0,2]$,
see Sullivan \cite{sull}; in fact, a $\delta_\ast$-conformal
measure exists for $\delta_\ast:=\inf\{\delta>0: \hbox{a conformal
measure $\mu_\delta$ exists}\}$.

Given $\delta>0$ and $z\in \C\sm\cup_{n \ge 1} f^n(0)$, let the
\emph{Poincar\'e series} be defined as
$$\Xi_\delta(z):=\sum_{n=0}^\infty\sum_{f^n(y)=z}\frac
1{|Df^n(y)|^\delta}.$$ As explained in \cite{GS}, a bounded
distortion argument implies that this number is independent of
$z\in\C\sm \omega(0)$. We let
$\delta_{Poin}:=\inf\{\delta>0:\Xi_\delta(z)<\infty \hbox{ for }
z\in \C\sm\omega(0)\}$.  For more information on the \poin series
see, for example, \cite{AL1,GS}.

We denote the Hausdorff dimension of a set $X$ by $HD(X)$.  The
Hausdorff dimension of a measure $\mu$ on $X$ is
$HD(\mu):=\inf\{HD(A):A\subset X, \mu(A)=1\}$. We define the
\emph{dynamical dimension} $DD(X):=\sup\{HD(\mu)\}$ where the
supremum is taken over all ergodic invariant measures of positive
entropy.  By \cite{PU}, $\delta_\ast=DD(\J)$, so clearly
$HD(\J)\ge \delta_\ast$. Also, by \cite{Bishop}, if the Lebesgue
measure of the Julia set is zero then $\delta_{Poin}\ge HD(\J)$.
Whenever $\delta_\ast=\delta_{Poin}$, $\mu_{\delta_\ast}$ is
called a \emph{geometric measure}.

We consider Feigenbaum and Fibonacci maps in this paper, which are
interesting examples in the measure theoretical context,
especially when the critical order $\ell$ is large. In this paper
a Feigenbaum map is an infinitely renormalisable map with periodic
combinatorics, i.e., with periodic renormalisation operator, see
\cite{MS}. In the Fibonacci case, the questions whether
$\delta_\ast < 2$, whether $\mbox{Area}(\J(f)) = 0$ and the
conservativity/dissipativity of $\delta$--conformal measure remain
unsolved. Avila and Lyubich \cite{AL1} show that in the Feigenbaum
case, $\mu_{\delta_\ast}$ is dissipative and that
$\mbox{Area}(\J(f)) = 0$ implies $\delta_\ast < 2$ (as well as
further results on a class of maps similar to Feigenbaum maps).

In the following theorem, proved in Section~\ref{sec:fib n feig},
we show that for Feigenbaum and Fibonacci-like (i.e., satisfying
\eqref{eq:fibo-like}  below) polynomials admit no acip $\mu \ll
\mu_\delta$ provided $\ell$ is sufficiently large.

\begin{theorem}\label{thm:fib no acip}
Suppose that $f:z \mapsto z^\ell+c_1$ is either a complex
Fibonacci-like map satisfying \eqref{eq:strictHofbauer} or a
Feigenbaum map.  If $\ell$ is sufficiently large, then $f$ cannot
have an acip w.r.t. any conformal measure, unless the acip is
supported on $\omega(0)$.
\end{theorem}

Note that $HD(\omega(0))\le 1$, so if
$\mu_{\delta_\ast}(\omega(0))>0$ then $HD(\J)=1$.  By
\cite{Zdparab}, this only occurs when $\J$ is a Jordan curve.  So
in particular, the maps in Theorem~\ref{thm:fib no acip} have no
acip w.r.t. $\mu_{\delta_\ast}$.

In fact, in \cite{AL1} it is proved that for quadratic Feigenbaum
maps, any conformal measure is dissipative (indeed this is claimed
for unicritical Feigenbaum maps of any critical order).
Lemma~\ref{lem:diss no acip} then implies that there cannot be an
acip.  However, we include the Feigenbaum case in our proofs for
interest. Further examples of maps having no acip are presented in
the following theorem.

\begin{theorem}
There exists a class of maps `long--branched maps' each of which
has no acip not supported on $\omega(0)$, for any even critical
order. \label{thm:lb no acip}
\end{theorem}

The paper is arranged as follows.  We start by outlining some
useful combinatorial facts for interval maps in
Section~\ref{sec:prelims}.  In Section~\ref{sec:no acip setup}, we
explain the general complex setup, and what conditions we require
to show there are no acips. In Section~\ref{sec:fib n feig}, we
show that these conditions are satisfied by certain Fibonacci and
Feigenbaum maps.  In Section~\ref{sec:long branches} we show that
these conditions are satisfied by some `long--branched' maps.
Section~\ref{sec:eq} applies the above results to the question of
(non--)existence of equilibrium states.  In the appendices we first
give some useful results on dissipativity/conservativity, and then
we supply some technical results on the kneading map.

Throughout we will write $x \approx y$ if there is a constant $K$
depending only on $f$ such that $\frac1K x \le y \le K x$.

{\bf Acknowledgements:} We would like to thank Juan
Rivera--Letelier, Feliks Przytycki and the referees for useful comments.

\section{Preliminaries}
\label{sec:prelims} Let $f(z) = z^\ell+c_1$ for $\ell \in 2\N$ and
$c_1 \in \R$, so the critical point is $0$ and
\[
f(e^{2\pi k i/\ell}\ \R) =\left\{ \begin{array}{ll}
[c_1,\infty) & \mbox{ if } k \mbox{ is even;}\\
(-\infty, c_1] & \mbox{ if } k \mbox{ is odd.}
\end{array} \right.
\]
Any maximal closed interval $J \subset \R$ such that $f^n|_J$ is
monotone is called a {\em branch} and if $0 \in
\partial J$, then the branch is called {\em central}; there are
two central branches one on either side of $0$, denoted $J_n$ and
$\hat J_n$. The numbers $1=S_0 < S_1 < S_2 < \dots$ such that if $n=S_k$ then
\[
D_n := f^n(J_n) \owns 0,
\]
are called the {\em cutting times} of $f$. It can be shown, see
\cite{Hof,BruinKnea}, that if there is no periodic attractor, then
there are infinitely many cutting times, and there is a function
$Q:\N \to \N \cup \{ 0 \}$, called the {\em kneading map}, such
that
\[
S_k - S_{k-1} = S_{Q(k)}\mbox{ for all } k \in \N \mbox{ (and
hence } Q(k) < k).
\]
As $D_{S_k} \owns 0$, there are points $\z_k < 0 < \hat \z_k$ such
that $f(\z_k) = f(\hat \z_k)$ and $f^{S_k}(\z_k) = 0$. Moreover,
there is no $z \in (\z_k,0)$ such that $f^n(z) =0$ for some $n \le
S_k$. For this reason, the points $\z_k$ and $\hat \z_k$ are
called {\em closest precritical points}.

\begin{lemma}\label{lem:basic props}
The following properties hold:
\begin{itemize}
\item[(a)] the central branches of $f^{S_k}$ are $J_{S_k} = [\z_{k-1},0]$ and
$\hat J_{S_k} = [0, \hat \z_{k-1}]$;
\item[(b)] the point $f^{S_k}(0)$ lies in $(\z_{Q(k+1)-1}, \z_{Q(k+1)})
\cup (\hat \z_{Q(k+1)}, \hat \z_{Q(k+1)-1})$.
\end{itemize}
\end{lemma}

\begin{proof} See \cite{BruinKnea}. \end{proof}

Simple examples of the kneading map are $Q(k) = k-1$, when $S_k =
2^k$ and the corresponding map is a Feigenbaum map; and $Q(k) =
\max\{ 0, k-2\}$, when $S_k$ are the Fibonacci numbers, and $f$ is
called the Fibonacci map. More generally, we call a map {\em
Fibonacci--like} if there exists $N$ such that
\begin{equation}\label{eq:fibo-like} k-Q(k) \le N \mbox{ for all
} k \ge 1.
\end{equation}
A Feigenbaum map is Fibonacci-like according to this definition,
but not all infinitely renormalisable maps are. In fact, (real)
renormalisability can be expressed in terms of the kneading map as
follows: There exists $k_0$ such that
\[
Q(k_0) = k_0-1 \mbox{ and } Q(k) \ge k_0-1 \mbox{ for all } k \ge
k_0.
\]
In this case, $f$ has a periodic interval of period $S_{k_0}$, see
\cite{BruinKnea}. The infinitely renormalisable maps that we are
interested in are those for which the renormalisation operator is
periodic. Hence, some iterate of the renormalisation operator will
fix the combinatorial structure, and this can be expressed by
periodicity in the kneading map: There exists $k_0$ and period $p$
such that
\begin{equation}\label{eq:feig-like}
Q(k_0) = k_0-1, Q(k) \ge k_0-1 \mbox{ and } Q(k+p) = Q(k)+p \mbox{
for all } k \ge k_0.
\end{equation}
Not every kneading map $Q$ corresponds to a real unimodal
polynomial. Hofbauer, in {\em e.g.}, \cite{Hof}, states the
following sufficient condition for a kneading map to be realised
by a polynomial.
\begin{equation} \{ Q(k+j) \}_{j \ge 1} \succeq
\{ Q(Q^2(k)+j) \}_{j \ge 1} \mbox{ for all } k \ge 1,
\label{eq:admiss}
\end{equation}
where $\succeq$ indicates lexicographical order and $Q^2(n)=Q(Q(n))$, and so on. A slightly stronger version of this
condition is as follows: there is $k_0$ such that
\begin{equation}\label{eq:strictHofbauer}
Q(k+1) > Q(Q^2(k)+1) \quad \mbox{ for all } k > k_0.
\end{equation}
We will use this condition, which excludes the existence of
``almost saddle node'' bifurcations, in Section~\ref{sec:fib n
feig} to simplify some of our results.

\section{Conditions which preclude acips}
\label{sec:no acip setup}

\subsection{The first return map to a wake}
\label{sec:kac}

We suppose throughout that $f:z \to z^\ell+c_1$ for $c_1\in \R$
has connected Julia set. We suppose further that $0$ is recurrent
and $0\in \J$ (if $0\notin \J$, then $f|_\J$ is hyperbolic and an
acip exists). This implies that the filled Julia set $K(f):=\{z\in
\C:f^n(z) \not\to \infty\}$ has no interior.
We let $r_{\theta}:=\{re^{i\theta}:r\in [0,\infty)\}$. Then for
$\theta \in T_\ell:= \left\{\frac \pi \ell,\frac{3\pi}\ell,
\ldots, \frac{(2\ell-1)\pi}\ell\right\}$ we have a connected
segment $r_\theta\cap\J$ mapping onto $\R^-\cap\J$. For example,
in the quadratic case, $i\R^+\cap\J$ and $i\R^-\cap\J$ map onto
$\R^-\cap\J$.

Let the Green function $G:\F \to \R$ be defined by $G(z) = \lim_{n
\to \infty} \frac{\log |f^n(z)|}{\ell^n}$. The
\emph{equipotentials} (i.e. level sets of the Green function) form
a foliation of the \emph{Fatou set} $\F$ consisting of nested
Jordan curves, see \cite{miln}. The orthogonal foliation is the
foliation of \emph{external rays}.  Let $B(\J)$ be the bounded
connected component of $\C\sm G^{-1}(\{ 1 \})$. Each external ray
$R$ has its \emph{external angle} $\alpha = \lim_{r \to \infty}
\arg(R \cap \{ |z| = r\})$, so we write $R=R_\alpha$. Given
external rays $R_\alpha, R_{\alpha'}$ landing at the same point
$z$, the corresponding \emph{wake} is the set of points in $B(\J)$
lying between $R_\alpha$ and $R_{\alpha'}$.

\begin{figure}[ht]
    \centering
    \psfragscanon
    \psfrag{label1}{{\Huge $c$}}
    \psfrag{label2}{{\Huge $0$}}
    \psfrag{labelW}{{\Large $W$} }
    \psfrag{labelWw}{{\Large $W_w$}}
    \scalebox{0.5}{\includegraphics{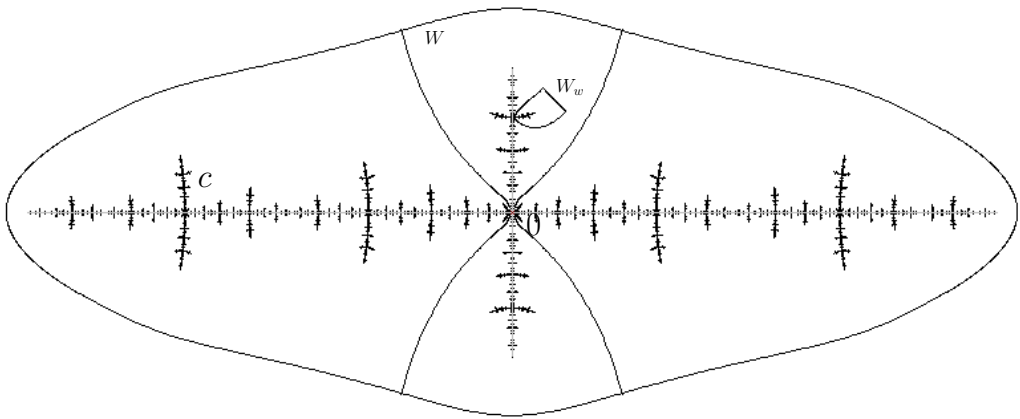}}
\caption{The neighbourhood $B(\J)$ of the Julia set with the wake
$W  = W^1 = W(\frac{\pi}{2})$. Here the map $f$ is quadratic, and
the first return map to $W$ maps the subwake $W_w$ onto $W$. The
points $0$ and $c$ are shown as well.} \label{fig:partition}
\end{figure}

Let $0\le\alpha<\frac{2\pi}\ell$ be such that $R_{\ell\alpha}, R_{\ell\tilde\alpha}$
are external rays landing at $c_1$ where $\tilde \alpha =
-\alpha$.  We pull these rays back by $f$ to get the external rays
$R_{\alpha_0}, R_{\tilde\alpha_0}, R_{\alpha_1},R_{\tilde\alpha_1}
,\ldots,R_{\alpha_{\ell-1}},R_{\tilde\alpha_{\ell-1}}$ where
$\alpha_k=  \alpha+\frac{2\pi k}\ell$ and $\tilde\alpha_k=
-\alpha+ \frac{2\pi k}\ell$.  We denote the wake corresponding to
$\tilde\alpha_k$ and $\alpha_k$ by $W^k$, see
Figure~\ref{fig:partition}. There are $\ell$ such wakes.  Each
$\theta\in T_{\ell}$ is of the form
$\theta=\frac{(2k+1)\pi}{\ell}$, so $r_\theta\cap\J\subset W^k$
and we write $W(\theta)=W^k$.

Fix $\theta\in T_\ell$ for the moment and consider the first
return map to a wake $W = W(\theta)$. We may index a domain of
this map by the point $w$, called \emph{root point}, that maps to
$0$.  We denote the set of root points by $R$.

The following proposition is a tool we use throughout the rest of
this paper.

\begin{proposition}
Let $W=W(\theta)$ be a wake where $\theta\in T_\ell$ and suppose
$\tau(z):=\min\{k\ge 1:f^k(z)\in W(\theta)\}$, wherever this is
defined.
If
\[
\sum_{w\in R}\tau(w)\mu_\delta(W_w)=\infty,
\]
then there is no acip with respect to $\mu_\delta$, unless the
acip is supported on $\omega(0)$. \label{prop:Kac no acip}
\end{proposition}

We will show that in fact the above sum is infinite when it is
taken over the set of domains whose root point $w$ lies on
$r_\theta$.

To prove this proposition, we first need the following result, the
proof of which is postponed to Appendix A.

\begin{lemma}
Suppose that $\mu$ is an acip w.r.t. $\mu_\delta$ and
$\mu(\omega(0))<1$. Then $\mu(W(\theta))>0$ for all $\theta\in
T_\ell$. \label{lem:cons}
\end{lemma}

\begin{lemma}[Koebe Distortion Lemma for complex maps] Assume
that $\disk$ is the open unit disk in $\C$ and $g\colon \disk\to
\C$ is univalent, then for each $z\in \disk$,
\[
\frac{1-|z|}{(1+|z|)^3} \le \frac{|Dg(z)|}{|Dg(0)|} \le
\frac{1+|z|}{(1-|z|)^3}.
\]
and
\[
\frac{1-|z|}{1+|z|} \le \frac{|z|\cdot |Dg(z)|}{|g(z)-g(0)|} \le
\frac{1+|z|} {1-|z|}.
\]
Thus, supposing that $z,w\in r\disk$ for some $0<r<1$, we have
\begin{equation}
\frac{|Dg(z)|}{|Dg(w)|} \le \frac{(1+r)^4}{(1-r)^4}.
\end{equation}
\label{lem_complex_koebe}
\end{lemma}

\begin{proof}
This lemma can be found in \cite{Pom}: the first statement is
formula (15) of Theorem 1.3, page 9, and the second follows by
substituting the Koebe transform $h(z) = ( g(\frac{z + w}{1+\bar w
z}) - g(w))/(1-|w|^2) Dg(w)$ in formula (14) and then taking $z =
-w$, cf. Exercise 3 on  page 13 in \cite{Pom}.
\end{proof}

\begin{proof}[Proof of Proposition~\ref{prop:Kac no acip}.]
Suppose that we do have an acip $\mu$ and $\mu(\omega(0))<1$.
\begin{claim} There is a lower bound on the density $b:=\inf_{x\in
W}\frac{d\mu}{d\mu_\delta}(x)>0$.\label{claim:bounded
density}
\end{claim}
\begin{proof}
First note that by Lemma~\ref{lem:cons} we have $\mu(W)>0$.  We
consider a domain $Q=W_w$ with $\mu(Q)>0$. There will be some
Koebe space for the map $f^n:Q \to W$ where $\tau|_Q=n$.  So there
exists $K$ \st $\left|\frac{Df^n(z)}{Df^n(v)}\right|< K$ for all
$z,v \in Q$. Fix this $v= v(Q)$ once and for all. We may assume
that the first return map $F_Q:\bigcup_i Q_i \to Q$ also has
distortion bounded by $K$. Let $\mu_{\delta,Q} :=
\frac{\mu_\delta|_Q}{\mu_{\delta}(Q)}$ and $\mu_Q :=
\frac{\mu|_Q}{\mu(Q)}$ be the normalised restrictions to $Q$ of
$\mu_\delta$ and $\mu$ respectively.  Note that since $\mu$ is an
invariant probability measure, Poincar\'e recurrence implies that
the set of points which return infinitely often to $Q$ has
positive $\mu$--measure, and thus positive $\mu_\delta$ measure.
Since the set of points which enter $Q$ infinitely often is an
invariant set with positive $\mu_\delta$ measure, the ergodicity
of $\mu_\delta$ given by \cite{Prado} means that this set has full
measure.

We can now apply the Folklore Theorem to $(F_Q, \mu_{\delta,Q})$, see for example \cite{MS},
and obtain an ergodic invariant measure $m_Q$. Then there exists
some $b'>0$ \st $\frac{dm_Q}{d\mu_{\delta,Q}}>b'$. But since $m_Q$
is ergodic and invariant and $\mu \ll \mu_\delta$, we have $m_Q =
\mu_Q$. Hence $\frac{d\mu_Q}{d\mu_{\delta,Q}}>b'$. Now consider a
domain $Q'= W_{w'}$ of $f^\tau$.  For any $U'\subset Q'$ there
exists a unique $U \subset Q$ \st $f^n(U)= U'$. By the
conformality of $\mu_\delta$ and the invariance of $\mu$,
\begin{align*} \mu_\delta(U') & = \int_{U}|Df^n|^\delta d\mu_\delta
\le K |Df^n(v)|^\delta \mu_\delta(U)  \le \frac{K}{b'}
|Df^n(v)|^\delta \mu(U) \frac{\mu_\delta(Q)}{\mu(Q)}\\
& \le \frac{K}{b'} |Df^n(v)|^\delta \mu(U')
\frac{\mu_\delta(Q)}{\mu(Q)}.
\end{align*}
Taking $b:=
\frac{b'}{K|Df^n(v)|^\delta}\frac{\mu(Q)}{\mu_\delta(Q)}$, the
proof of the claim is finished.
\end{proof}
Since $\mu$ is invariant, Kac's Lemma implies
$$
1 \ge \sum \tau(w) \mu(W_w) \ge b\sum \tau(w) \mu_\delta(W_w),
$$
contradicting the assumptions of our lemma.  So there is no acip.
\end{proof}

\subsection{Wake boundaries}
\label{sec:wakes}

To apply Proposition~\ref{prop:Kac no acip} we will need to
estimate $\mu_\delta(W_w)$ in terms of $|Df^{\tau(w)}(w)|$.  We do
this by establishing distortion bounds for the first return map to
a suitably chosen wake $W$ for some domains $W_w$ and then
applying Lemma~\ref{lem_complex_koebe}.

In order to find the necessary Koebe space, we make use of the
fact that any forward iterate of the critical point for our maps
$f:z \mapsto z^\ell+c_1$ must lie in the real line and that the
wake $W$ is separated from certain iterates of the critical point
as explained below. In Section~\ref{sec:fib n feig} we need $W$ to
be separated from $\R\sm\{0\}$, and in Section~\ref{sec:long
branches} for a neighbourhood $U$ of 0 we need $W\setminus U$ to
be a bounded distance away from $\R$.

\begin{lemma}
Suppose that $f(z) = z^\ell+c_1$ for $c_1\in \R$ and $\ell \ge 6$.
Then
\[
W^k \subset \left\{re^{i\theta}:r\ge 0,\ \theta \in \left(\frac
{2\pi} \ell, \frac{(\ell-1)\pi} \ell\right)\right\}
\]
whenever $2\le k\le \frac{\ell-1}2$. \label{lem:good angles}
\end{lemma}

\begin{proof}
Given $W^k$, let $\W^k$ denote the wake constructed as $W^k$ was,
but not bounded by an equipotential.  So $W^k\subset \W^k$ and on
the Riemann sphere $\W^k\cap\W^j =\{\infty\}$ for $j\neq k$.  We
prove the lemma for $\W^k$ which implies that it holds for $W^k$
too.

We fix $r>0$ and let $C_r:=\{re^{i\theta}:\theta\in [0,2\pi)\}$.
Let $\alpha_k:=\sup\{\theta\in[0,2\pi):re^{i\theta}\in \bd \W^k\}$
and $\beta_k:=\inf\{\theta\in[0,2\pi):re^{i\theta}\in \bd \W^k\}$.
Notice that $\beta_k=\beta_{k-1}+\frac{2\pi}\ell$ and
$\alpha_k=\alpha_{k-1}+\frac{2\pi}\ell$.
\begin{claim}
There is no $r>0$ such that $re^{i\xi_{k+1}} \in \W^{k+1}$ with $\xi_{k+1}< \beta_k$.
\end{claim}
Suppose that the claim is false.  Then there exists $r>0$ such that $re^{i\xi_{1}} \in \W^{1}\cap C_r$ with $\beta_1 \le\xi_1< \beta_0$. Since $\beta_0<0<\beta_1$ (no ray may cross $\R$), we have a
contradiction, proving the claim.

Since $\alpha_k, \beta_k>\frac{2\pi}\ell + \beta_1$ for $k\ge 2$,
the claim implies
$$\W^k \subset \left\{re^{i\xi}:r\ge 0,\
\xi \in \left(\frac{2\pi} \ell,
\frac{(\ell-1)\pi}\ell\right)\right\}$$ whenever $\alpha_k \le
\frac{(\ell-1)\pi}\ell$. Since $\alpha_k \le \frac{2\pi k}\ell$,
the lemma is proved.
\end{proof}

\begin{lemma}
Suppose that $U$ is a neighbourhood of $0$.   Then for any ray
$\gamma$ landing at 0 there exists $\eps=\eps(\gamma)$ \st
$dist(\gamma\sm U, \J\cap\R)>\eps$, where $dist$ denotes distance in the Euclidean metric. \label{lem:LCJulia set}
\end{lemma}

\begin{proof}
First notice that $\J$ is symmetric with respect to both $\R$ and
$i\R$.  So if there is an external ray landing at $0$ in one
quadrant of $\C$ then there must be one in all other quadrants of
$\C$.  Therefore, close to $\J$, for example inside the
equipotential $\{ G(z) = \delta\}$ for small enough $\delta>0$, a
ray landing at 0 must start in a given quadrant of $\{ G(z) =
\delta\}$ and remain in that quadrant until it lands at 0.  Of
course, rays cannot intersect $\J$.

Now suppose that the lemma is false. Then there is some sequence
of points $y_n\in (\J\cap\R)\sm U$ \st $d(y_n, \gamma) \to 0$ as
$n \to \infty$. But this implies that for rays to land at points
in $\J\cap\R\cap U$ they must intersect $\J\cap\R$ which is not
possible. Since $\J$ is locally connected, see \cite{LvSlocconn},
every point of $\J$ has a ray landing by Carath\'eodory's Theorem.
So we have a contradiction.
\end{proof}

\section{Fibonacci-like and Feigenbaum maps}
\label{sec:fib n feig}

This section is devoted to the proof of Theorem~\ref{thm:fib no
acip}.

\begin{proposition}\label{prop:scaling}
Assume that $Q$ is an admissible kneading map satisfying either
\eqref{eq:fibo-like} and \eqref{eq:strictHofbauer}
(Fibonacci-like), or \eqref{eq:feig-like} (Feigenbaum). Then for
each $\ell \in 2\N$, there is a unique $c_1$ such that $f(z)
=z^\ell + c_1$ has kneading map $Q$. If $\ell > 2$, there exists
$\lambda=\lambda(\ell)<1$ such that
\[
\frac{|f^{S_{k+1}}(0) - 0|}{|f^{S_k}(0) - 0|} \to \lambda, \mbox{
as } k \to \infty  \hbox{.  Moreover, } \lambda \to 1 \mbox{ as }
\ell \to \infty.
\]
\end{proposition}

\begin{proof} The existence of the parameters $c_1 = c_1(\ell)$
follows from the admissibility of the kneading map $Q$. Uniqueness
of $c_1(\ell)$ follows from rigidity results due to \cite{L,GS}
for $\ell = 2$ and \cite{KSS} for $\ell
> 2$. Finally the scaling properties were proved in \cite{BKNS}
for the Fibonacci map, generalised to the Fibonacci--like case in
\cite[Proposition 10.6]{BruinAttr}.  The Feigenbaum  case follows
from \cite{LevSwi}.
\end{proof}

Let $S_n$ denote the $n$--th cutting time (for example, the
$n$--th Fibonacci number for the Fibonacci map).  We will need the
following technical lemma.

\begin{lemma}
Suppose that the critical point of $f$ has order $\ell$. Assume
that the kneading map $Q$ of $f$ either satisfies
\eqref{eq:fibo-like} and \eqref{eq:strictHofbauer}
(Fibonacci-like) or \eqref{eq:feig-like}  (Feigenbaum).
Then there exists $B=B(\ell)>0$ \st given $\theta \in T_\ell$ \st
$W(\theta) \subset\left\{re^{i\xi}:r\ge 0,\ \xi \in
\left(\frac{2\pi} \ell, \frac{(\ell-1)\pi}\ell\right)\right\}$,
for each $n$ there is $w_n \in r_\theta$ such that
$f^{S_n}(w_n)=0$ and $|Df^{S_{n}}(w_{n})|<B$. \label{lem:bounded
derivative}
\end{lemma}

\begin{proof}
Let $(a,b)$ denote the open interval with endpoints $a$ and $b$,
even if $b <a$.

Recall that $\z_n<0<\hat \z_n$ are \st $f^{S_n}(\z_n)=
f^{S_n}(\hat \z_n) = c$, and $f^{S_n}:(\z_n,0) \to (0,
f^{S_n}(0))$, $f^{S_n}: (0,\hat \z_n) \to (0, f^{S_n}(0))$ are
monotone.  Let $w_n\in r_\theta$ be the point which has
$f(w_n)=f(\z_n)$.

We give bounds for derivatives $|Df^{S_n}(c_1)|$,
$|Df^{S_n}(\z_n)|$ and $|Df^{S_n-1}(c_1)|$. By
Lemma~\ref{lem:monotone} in Appendix B there is some interval $V$
containing $c_1, f(w_n), f(\z_{n}) $ \st the map $f^{S_n-1}: V \to
(f^{S_{Q(n)}}(0), f^{S_{Q^2(n)}}(0))$ is monotone, as in
Figure~\ref{fig:derivatives}.

\begin{figure}[ht]
\unitlength=11mm
\begin{picture}(10,5)(0,0)
\thinlines \put(0,4){\line(1,0){10}}
\put(4.5,4.3){\vector(-1,0){4.3}}\put(4,4.3){\vector(1,0){5.5}}
\put(4.5,4.5){\small $V$} \put(0.2,3.9){\line(0,1){0.2}}
\put(5,3.9){\line(0,1){0.2}} \put(4.8,3.5){\tiny $c_1$}
\put(6.5,3.9){\line(0,1){0.2}} \put(6.3,3.5){\tiny $f(\z_{n+1})$}
\put(8,3.9){\line(0,1){0.2}} \put(7.8,3.5){\tiny $f(\z_{n})$}
\put(9.5,3.9){\line(0,1){0.2}} \put(9.2,3.5){\tiny $f(\z_{n-1})$}
\put(0,3.5){\vector(0,-1){1.4}} \put(-1.3,2.5){\small
$f^{S_{n}-1}$}

\put(0,1.5){\line(1,0){10}} \put(0.2,1.4){\line(0,1){0.2}}
\put(-0.4,1){\tiny $f^{S_{Q^2(n)}}(0)$}
\put(5,1.4){\line(0,1){0.2}} \put(4.6,1){\tiny $f^{S_{n}}(0)$}
\put(6.5,1.4){\line(0,1){0.2}} \put(6.1,1){\tiny $\z_{Q(n+1)}$}
\put(8,1.4){\line(0,1){0.2}} \put(8,1){\tiny $0$}
\put(9.5,1.4){\line(0,1){0.2}} \put(9,1){\tiny $f^{S_{Q(n)}}(0)$}

\end{picture}
\caption{Koebe space for $f^{S_n-1}$ to show that
$|Df^{S_{n}-1}(f(\z_n))| \approx |Df^{S_{n}-1}(c_1)|$.}
\label{fig:derivatives}
\end{figure}
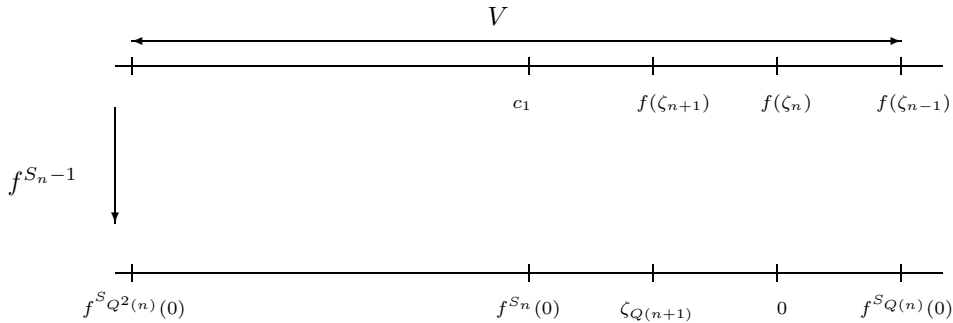

Moreover, $(0, f^{S_{n}}(0))$ is well inside $(f^{S_{Q^2(n)}}(0),
f^{S_{Q(n)}}(0))$; indeed by Proposition~\ref{prop:scaling} and
assumption \eqref{eq:strictHofbauer} it is a $(1-\lambda)$--scaled
neighbourhood of $(0,f^{S_{n}}(0))$ (note that $\lambda$ depends
on $\ell$). Thus we can use the Koebe lemma to get
\[
|Df^{S_{n}-1}(c_1)| \approx
\frac{|f^{S_{n}}(0)-0|}{|f(\z_{n})-c_1|} \approx
\frac{|f^{S_{n}}(0)-0|}{|\z_{n}-0|^\ell},
\]
recall that we write $x \approx y$ if there is a constant $K$
depending only on $f$ such that $\frac1K x \le y \le K x$.
Hence
\[
|Df^{S_{n}}(c_1)| \approx
\frac{|f^{S_{n}}(0)-0|}{|\z_{n}-0|^\ell}|f^{S_{n}}(0)-0|^{\ell-1}
= \left(\frac{|f^{S_{n}}(0)-0|}{|\z_{n}-0|}\right)^\ell.
\]
Since, by Lemma~\ref{lem:basic props}(b), $f^{S_n}(0) \in
(\z_{Q(n+1)-1}, \z_{Q(n+1)}) \subset (\z_{Q(n+1)-1}, \z_{n})$, if
$Q$ satisfies \eqref{eq:fibo-like} then the right hand side above
is bounded above by
$\left(\frac{|\z_{n-N-1}-0|}{|\z_{n}-0|}\right)^\ell =O(1)$.
Similarly if $Q$ satisfies \eqref{eq:feig-like} then letting
$M:=Q(k_0)-k_0$, the right hand side above is bounded above by
$\left(\frac{|\z_{n-M-1}-0|}{|\z_{n}-0|}\right)^\ell =O(1)$ for
all $n\ge k_0$.

So in either case, for any large $n$ we have
\begin{align*} |Df^{S_{n}}(\z_{n})| & =
|Df^{S_{n}-1}(f(\z_{n}))Df(\z_{n})| \approx |Df^{S_{n}-1}(c_1)Df(\z_{n})|\\
& \approx
\frac{|f^{S_{n}}(0)-0|}{|\z_{n}-0|^\ell}|\z_{n}-0|^{\ell-1} = O(1)
\end{align*}
as required.
\end{proof}

\begin{proof}[Proof of Theorem~\ref{thm:fib no acip}]
Let $Pre$ be the set of points constructed in Lemma~\ref{lem:bounded
derivative}.  By Lemma~\ref{lem:good angles}, for $\ell \ge 6$ we
can choose $\theta$ so that the wake $W=W(\theta)$ described in
Section~\ref{sec:wakes} has $r_\theta\cap\J \subset W$ and
$W\subset\left\{re^{i\xi}:r\ge 0,\ \xi \in \left(\frac{2\pi} \ell,
\frac{(\ell-1)\pi}\ell\right)\right\}$. We now wish to estimate
$\mu_\delta(W_{w_n})$ using the derivative $Df^{S_n}(w_n)$ and
distortion properties of $f^{S_n}: W_{w_n} \to W$.

As in Figure~\ref{fig:derivatives}, there is some interval $V \ni
c_1, f(w_n)$ \st $f^{S_n-1}:V \to (f^{S_{Q^2(n) }}(0), f^{S_{Q(n)}}(0))$
is monotone.  Since $|f^{S_n}(0)-0| \ge C\lambda^n$ we have some
Koebe space for $f^{S_{n}}:W_{w_n} \to W$, which we use below.

We can trace out Koebe space for $f^{S_n}:W_{w_n} \to W$ using an
enlarged wake $W'$ whose boundary is parallel to this boundary
curve of $W$. By Proposition~\ref{prop:scaling} and since we chose
$W\subset\left\{re^{i\xi}:r\ge 0,\ \xi \in \left(\frac{2\pi} \ell,
\frac{(\ell-1)\pi}\ell\right)\right\}$, the Koebe space is greater
than $C\lambda^{n}$ where $C>0$ depends only on $\ell$. According
to Lemma~\ref{lem_complex_koebe}
applied to $f^{-S_n}:W' \to W'_{w_n}$, we have
$\mu_\delta(W_{w_n}) \ge K|Df^{S_n}(w_n)|^{-\delta}\lambda^{4n}$.
To compare with Proposition~\ref{prop:Kac no acip} we sum over the
points $w_{n}\in Pre$ obtained in Lemma~\ref{lem:bounded derivative}
to get
\begin{equation}\label{eq:upper} \sum_{w_{n}\in Pre} \tau(w)\mu_\delta(W_{w_{n}}) \ge K\sum_{w_n\in Pre} S_n
|Df^{S_n}(w_n)|^{-\delta}\lambda^{4n}>KB^{-\delta}\sum_{w_n\in Pre}
S_n\lambda^{4n}.
\end{equation}
For Fibonacci--like maps, $S_n$ increases exponentially, say $S_n
\asymp \sigma^n$ for some $\sigma > 1$ depending only on the
combinatorial type of $f$. For example, if $f$ is the Fibonacci
map, $S_n \asymp \gamma^n$ where $\gamma = (1+\sqrt{5})/2$ is the
golden ratio. If $\sigma \lambda^{4}>1$, then the proof is
completed by Proposition~\ref{prop:Kac no acip}. Recalling, from
Proposition~\ref{prop:scaling}, that $\lambda \to 1$ as $\ell \to
\infty$, we can increase the order of the critical point to reach
this conclusion. For Feigenbaum maps with periodic combinatorics
period $p$ we note that $S_n\asymp p^n$, so we may apply the same
argument as for the Fibonacci--like case.
\end{proof}

\section{Long--branched maps}
\label{sec:long branches}

A real map $f$ is called {\em long--branched} if there is $\kappa >
0$ such that for all $n$ and all maximal intervals $J \subset \R$
such that $f^n|_J$ is monotone, $|f^n(J)| > \kappa$. In
\cite{BrLongb}, these maps were studied in connection to cascade
of saddle-node bifurcations, and it was shown that under
appropriate combinatorial conditions, no invariant probability
measure can be absolutely continuous with respect to
one-dimensional Lebesgue measure. In this section, we prove a
similar non-existence result for probability measures that are
absolutely continuous with respect to the $\delta$-conformal
measure on the Julia set. As before, let $(S_k)_{k \ge 0}$ be the
cutting times, and let $\z_k$ and $\hat \z_k$ stand for the
closest precritical points.

\begin{theorem}\label{thm:lb1}
If $f(z) = z^{\ell} + c$, $c \in \R$, $\ell \in 2\N$, is
long--branched and
\begin{equation}\label{eq:lb1}
\sum_k S_k |\z_k - \z_{k+1}|^{-\delta} = \infty,
\end{equation}
then $f$ cannot have an acip with respect to any
$\delta$--conformal measure, unless the acip is supported on
$\omega(0)$.
\end{theorem}

\begin{remark}
If $f$ is long--branched, then it is straightforward  to show that
the kneading map $Q$ is bounded, say $Q(k) \le B$, (in fact, this
is equivalent to long--branchedness) and that $k-1 \le S_k \le S_B
k$ for all $k \ge 0$, see \cite{BrLongb}.
\end{remark}

\begin{lemma}
If $f$ is long--branched, then there exists $\kappa' > 0$ such
that $\kappa' \le |f^{S_k}(0) - \z_{Q(k+1)}| \le | f^{S_k}(0)-0|$
for all $k \ge 0$. \label{lem:longbr}
\end{lemma}

\begin{proof}
Suppose by contradiction that there is a sequence $\{ k_i \}_i$
such that  $|f^{S_{k_i}}(0) - \z_{Q(k_i+1)}| \to 0$. Since
$Q(k_i+1) \le B$, it follows that
\[
|f^{1+S_{Q(k_i+1)}} \circ f^{S_{k_i}}(0) -
f^{1+S_{Q(k_i+1)}}(\z_{Q(k_i+1)})| = |f^{ 1+S_{ k_{i+1} } }(0) -
f(0)| \to 0
\]
as well. But this contradicts that $f$ is long--branched.
\end{proof}

\begin{proof}[Proof of Theorem~\ref{thm:lb1}]
Fix some odd $m < \ell$. Let $r > 1$ be minimal such that $f(0) <
f^r(0) < 0$. Then $f^r: e^{\frac{\pi m i}\ell}\ \R \to [f^r(0),
\infty)$ is an onto $2$-to-$1$ map with a single  branch point at
$0$. Hence for each $\z_k \in (f^r(0), 0)$, there exists $w_k \in
e^{\frac{\pi m i}\ell}\R$ such that $f^r(w_k) = \z_k$, and the
corresponding domains $W_{w_k}$ have $f^r(W_{w_k}) = W_{\z_k}$,
see Section~\ref{sec:wakes}.

The map $f^{S_k}:(\z_{k-1}, 0) \to (f^{S_{Q(k)}}(0), f^{S_k}(0))$
is monotone onto, and by Lemma~\ref{lem:longbr}, its image $(f^{
S_{Q(k)} }(0), f^{S_k}(0))$ contains a $\kappa'$-neighbourhood of
$f^{S_k}([\z_k, \z_{k+1}]) = [0,\z_{Q(k+1)}]$. As $f$ has negative
Schwarzian derivative, the distortion of $f^{S_k}|_{[\z_k,
\z_{k+1}]}$ is bounded, say by $K = K(\kappa')$, independently of
$k$. We find
\[
|Df^{S_k}(\z_k)| \le K \  \frac{|\z_{Q(k+1)} -
0|}{|\z_k-\z_{k+1}|}.
\]
Since $|w_k-0| \approx |\z_k - f^r(0)| > 0$ for large $k$, we get
that $Df^r(w_k)$ is bounded and bounded away from $0$, uniformly
in $k$. Therefore, by Lemma~\ref{lem:longbr},
\[
|Df^{r+S_k}(w_k)|^{-\delta} \approx |\z_k-\z_{k+1}|^{\delta}.
\]
Recall that by Lemma~\ref{lem:basic props2}(a) $f^{S_k}$ maps
$(\z_{k-1}, 0)$ diffeomorphically onto
$(f^{S_k}(0),f^{S_{Q(k)}}(0))$.  Hence, by Lemma~\ref{lem:LCJulia
set}, there exists a neighbourhood $W'$ of $W$, with $\partial W'$
intersecting $\R$ at $f^{S_k}(0)$ and $f^{S_{Q(k)}}(0)$, and
neighbourhoods $W'_{w_k} \supset W_{w_k}$ \st $f^{r+S_k}:W_{w_k}'
\to W'$ is univalent. By Lemma~\ref{lem:longbr}, the Koebe space
for this map is of order $\kappa'$, whence
\[
\mu_\delta(W_{w^{\pm}_k}) \approx  |\z_k-\z_{k+1}|^{\delta}.
\]
Now since
\[
\sum_k \tau(w_k) \ \mu(W_{w_k}) \approx \sum_k (r+S_k)
|\z_k-\z_{k+1}|^{\delta} = \infty,
\]
the theorem follows by Proposition~\ref{prop:Kac no acip}.
\end{proof}

To show that there are indeed maps satisfying the conditions of
Theorem~\ref{thm:lb1}, we are going to specify combinatorial
conditions that imply condition \eqref{eq:lb1}. We suppose that
$f^N:J \to f^N(J)$ is the central branch and $(f^{kN}(0))_{k =
0}^d$ is a monotone sequence of points in $J$. When $d$ is large
(so $f^N:J \to f^N(J)$ is ``very close to the diagonal''), we call
this an \emph{almost tangency}. The largest $d \in \N$ such that
$f^{dN}(0) \in J$ is called the {\em depth } of the
almost-tangency. We speak of a {\em cascade of almost saddle-node
bifurcations} if there is a sequence $(N_j)_{j \in \N}$ of
iterates such that the graphs of the central branches of $f^{N_j}$
are close to tangency for all $j$. Let $(d_j)_{j \in \N}$ be the
corresponding depths. Note that given $N_j$, one can adjust the
parameter $c$ such that $d_j$ can be arbitrarily large, see
\cite{BrLongb}.

\begin{proposition}
If $f$ is long--branched, and has an almost saddle node cascade
satisfying
\[
\sum_j \frac{d_j}{N_j L^{2N_j/(\ell-1)} } = \infty,
\]
where $L := \sup\left\{|Df(x)| \ : \ x \in [f(0), f^2(0)] \right\}
\le 2\ell$, then condition \eqref{eq:lb1} holds for $\delta\le 2$.
\end{proposition}

\begin{proof}
Since the central branch of $f^{N_j}$ near $0$ consists of $x
\mapsto x^{\ell}$ composed with a diffeomorphism of bounded
distortion uniformly in $j$, we can write
\[
f^{N_j}(x) \approx g(x) := \alpha x^{\ell} + \beta,
\]
for $x$ close to $0$.

Without loss of generality we can assume that $f^{N_j}$ has a
local maximum at $0$. Let $[z,0]$ be the domain of this central
branch, so $f^n(z) = 0$ for some $n < N_j$. Let $a\in \N$ be such
that $f^{N_j}(0) \in (\z_{a-1}, \z_a)$ and $f^{N_j}(\z_a) = \z_b$
is a closest precritical point too, so $N_j = S_a - S_b$, and
$f^{2N_j}(0) \in (\z_{b-1}, \z_b)$. Repeating this iteration, we
have $f^{d_jN_j}(\z_a) = \z_k$ and $f^{N_j}(\z_k)$ is still a
closest precritical point, but it lies outside $[z,0]$, so $S_k
\le 2N_j$. It follows that $d_j N_j < S_a < (d_j+2)N_j$.

Let $x_0$ be such that $Dg(x_0) = 1$, i.e., $x_0 = (\alpha
\ell)^{-1/(\ell-1)}$. Because $g$ has no fixed point,
\[
x_0 < g(x_0) = \alpha \left(\alpha \ell\right)^{-\ell/(\ell-1)} +
\beta = \alpha^{-1/(\ell-1)} \ell^{-\ell/(\ell-1)} + \beta.
\]
This gives
\[
\beta > \alpha^{-1/(\ell-1)} \left( \ell^{-1/(\ell-1)} -
\ell^{-\ell/(\ell-1)} \right) \ge \frac14 \alpha^{-1/(\ell-1)},
\]
and using $\alpha\le L^{N_j}$,
\[
g^2(0) - g(0) = \alpha \beta^\ell + \beta - \beta \ge 4^{-\ell}
\alpha^{-1/(\ell-1)} \ge 4^{-\ell} L^{-N_j/(\ell-1)}.
\]
There are at most $N_j$ integers $t$ such that $\z_t \in
(f^{N_j}(0), f^{2N_j}(0))$, so for at least one such $t$,
\[
|\z_t - \z_{t+1}| \ge \frac{4^{-\ell}}{N_j-1} L^{-N_j/(\ell-1)}.
\]
For this $t$, $S_t \ge (d_j-1) N_j$. Therefore we find
\[
\sum_k S_k |\z_k-\z_{k+1}|^\delta \ge \sum_j (d_j-1) N_j \left(
\frac{ 4^{-\ell} }{N_j-1} L^{-N_j/(\ell-1)}\right)^{\delta}
> 17^{-\ell} \sum_j d_j L^{ -2N_j/(\ell-1) }/N_j.
\]
By the assumption of the proposition, this sum diverges (and hence
no acip exists).
\end{proof}

\section{Equilibrium States}
\label{sec:eq}

We show the implications of the results in the previous section to
the question of the existence of equilibrium states for the
potential $-\delta\log |Df|$, i.e., whether or not there
exists a measure which achieves the supremum (called
\emph{pressure})
$$
P := \sup_\mu \left( h_\mu-\delta \int_{\J}\log|Df|~d\mu
\right)
$$
taken over set of ergodic invariant Borel probability measures.
This question is answered by Proposition~\ref{prop:no eq}. We need
the following result from \cite{Ledcomp} to prove our proposition.

\begin{theorem}
Suppose that $f$ has a $\delta_\ast$--conformal measure and $\mu$
is an invariant measure with $\lambda(\mu)>0$. Then
$$h_\mu=\delta_\ast\int\log|Df|~d\mu \quad  \hbox{ if and only if }
\quad\mu\ll\mu_{\delta_\ast}.$$ \label{thm:led}
\end{theorem}

\begin{proposition}
Suppose that $h_\mu=\delta_\ast\int\log|Df|~d\mu$ for a
unicritical map $z\mapsto z^\ell+c$ for $c\in \R$ and $\ell\in
2\N$.  If $\lambda(\mu)>0$ then $\mu\ll\mu_{\delta_\ast}$.
Otherwise $\mu$ is in the convex hull of weak accumulation points
of $\frac 1{n}\sum_{k=0}^{n-1}\delta_{f^k(0)}$. \label{prop:no eq}
\end{proposition}

\begin{proof}
By Theorem A of \cite{Przlyap}, for any invariant measure $\mu$
supported on $\J$, $\lambda(\mu)\ge 0$. First suppose that
$\lambda(\mu)>0$. According to for example \cite{DU} or Theorem
A2.9 of \cite{Przcon}, the pressure is zero. So we must have
$h_\mu=\delta_\ast\int\log|Df|~d\mu$. If $\lambda(\mu)>0$ then
Theorem~\ref{thm:led} implies $\mu\ll\mu_{\delta_\ast}$.

If $\lambda(\mu)=0$ and $\mu(\J \setminus \R) = 0$, then the
proposition is completed by the result of Hofbauer and Keller
\cite{HK} for interval maps which states that $\mu$ belongs to the convex hull of
the weak accumulation points of $\frac1n \sum_{k=0}^{n-1}
\delta_{f^k(0)}$. Otherwise, we can choose a set $A$ \st
$\dist(A,\R)>0$ and $\mu(A)>0$.

Given $\rho  > 0$ and $M > 0$ we say that $z$ \emph{reaches large
scale} at time $i$ if there are neighbourhoods $\C \supset V_0
\supset V_1 \owns z$ such that $f^n:V_0 \to f^n(V_0)$ is
univalent, $f^i(V_1)$ contains a round ball of radius $\rho$
(measured in Euclidean distance) and the modulus of
$\modulus(V_0,V_1) > M$.

Since $\mu$ is ergodic and invariant, this means that $\mu$--a.e.
$z$ visits $A$ with positive frequency, and hence goes to large
scale (with bounded distortion) with positive frequency. The paper
\cite{BT} then implies that we can `lift' $\mu$, which implies
that $\mu$ has positive Lyapunov exponent, which is a
contradiction.
\end{proof}

\begin{corollary}
The maps considered in Theorems~\ref{thm:fib no acip} and
\ref{thm:lb no acip} may only have an equilibrium state in the
convex hull of weak accumulation points of $\frac
1{n}\sum_{k=0}^{n-1}\delta_{f^k(0)}$. \label{cor:no eq}
\end{corollary}

For Fibonacci--like maps exactly one equilibrium state exists and
is supported on the minimal Cantor set $\omega(0)$, see
\cite{BruinKnea} and \cite{BrUE} for the unique ergodicity of
$f|_{\omega(0)}$.

\begin{proof}
For $\delta =\delta_\ast$, this follows directly from the
non--existence of acips, Theorem~\ref{thm:led} and
Proposition~\ref{prop:no eq}.  In the case $\delta>\delta_\ast$ we
have the following argument.  Suppose that $\mu$ is an an
equilibrium state.  First suppose that $\mu(\R)>0$.  By ergodicity
$\mu(\R)=1$.  Then the argument of \cite{BruinKnea} implies that
$\mu$ is in the convex hull of weak accumulation points of $\frac
1{n}\sum_{k=0}^{n-1}\delta_{f^k(0)}$.

Now suppose that $\mu(\R)=0$.  We can use the argument in the
proof of Proposition~\ref{prop:no eq} to show that
$\lambda(\mu)>0$.  Next we use a result of \cite{Ledcomp,
mane}, see also \cite[Theorem 9.4.1]{PU}, which states that
$HD(\mu)=\frac{h_\mu}{\lambda(\mu)}$.  Since by \cite{DU},
$DD(\J)=\delta_\ast$, we have $h_\mu-\delta\int\log|Df|~d\mu<0$.
However, \cite{BrK} implies that the pressure of $-\delta\log|Df|$
is greater than or equal to zero: a contradiction.
\end{proof}

\begin{corollary}
There exists a quadratic long--branched map with no equilibrium
state in its Julia set for the potential $z\mapsto
-\delta\log|Df(z)|$ for each $\delta \ge \delta_\ast$.
\label{cor:no eq for lb}
\end{corollary}

\begin{proof}
From the proof of Corollary~\ref{cor:no eq} we only need to
exclude the case that the equilibrium state is supported on $\R$.

In \cite{HK} an example of a real quadratic map was presented such
that $\frac1{n}\sum_{k=0}^{n-1}\delta_{f^k(x)}$ converges to the
Dirac measure at the (repelling) fixed point $p \in [f(0), 0]$ for
Lebesgue--a.e. $x \in [f(0), f^2(0)]$. In \cite[Example 5.4]{BrK}
this example was put in the context of long--branched maps, and it
was shown that it did not have an equilibrium state. The same
example works in the complex case. Corollary~\ref{cor:no eq}
implies that the only equilibrium state is the Dirac measure at
$\{ p \}$, so its free energy is $-\log |f'(p)| < 0$. On the other
hand, there is a sequence $\{ p_i\}_{i \ge 1}$ of periodic points
whose periods $\mbox{per}(p_i) \to \infty$ but whose multipliers
$|Df^{\mbox{\tiny per}(p_i)}(p_i)|$ remain bounded. So the
pressure is $P \ge 0$. This shows that no equilibrium state
exists.
\end{proof}

\section*{Appendix A} \label{sec:dissipative}

We follow \cite{AL1} in proving some facts about
dissipativity/conservativity.  While the ideas presented here are
very similar to those in \cite{AL1}, we do not make any assumption
on the combinatorics of our maps.  Notice that we restrict
ourselves to unicritical polynomials for simplicity, but all the
results below extend to general complex polynomials.  Recall that
$(X,T,\mu)$ is dissipative if there exists $A\subset X$ \st
$\mu(A)>0$ and $\mu(T^{-n}A\cap A)=0$ for all $n\ge 1$.

\begin{lemma}
Suppose that $f$ is a  unicritical polynomial with $\omega(0) \neq
\J$. Then the following are equivalent:
\begin{itemize}
\item[(a)] $\mu_\delta$ is dissipative;
\item[(b)] $\Xi_\delta(z)$ is convergent for all $z\notin\omega(0)$;
\item[(c)] $f^n(z)\to \omega(0)$ for $\mu_\delta$--a.e. $z\in\J$.
\end{itemize}
Moreover, if $\mu_\delta$ is conservative then almost all orbits
are dense in $\J$.\label{lem:diss}
\end{lemma}

\begin{proof}
Notice that $\omega(0)\neq \J$ implies that
$\mu_\delta(\omega(0))<1$.

(a)$\Rightarrow$ (b):  Let $Y$ be a wandering set of positive
measure.  Choose $r>0$ \st $B_{2r}(z)\cap \omega(0) =\es$. Let
$D=B_r(z)$. Since there exists $N\ge 0$ \st $f^N(D) \supset \J$,
there exists $Y'\subset D$ \st $f^N(Y') \subset Y$ and
$\mu_\delta(Y')>0$. By the Koebe Lemma, and the fact that
preimages of $Y'$ are pairwise disjoint, we have
$$\Xi_\delta(z) \approx \sum_{n\ge 0}\mu_\delta(f^{-n}(Y'))<\infty$$
as required.

(b) $\Rightarrow$ (c): Take any $z\in \J$ and $r>0$ \st
$B_{2r}\cap \omega(0) =\es$.  Letting $D=B_r(z)$, we have
$$
\sum_{n\ge 0}\mu_\delta(f^{-n}(D))\approx \Xi_\delta(z)<\infty.
$$
By the \bcls, for $\mu_\delta$--a.e. $y\in \J$, the orbit of $y$
visits $D$ only finitely often.  Since this is true of any open
set bounded away from $\omega(0)$ this implies (c).

(c) $\Rightarrow$ (a): Let $U$ be a neighbourhood of $\omega(0)$
with $0<\mu_\delta(U)<1$.  Define $Y_n:=\{z\in \C:f^j(z)\in U
\hbox{ for all } j\ge n\}$.  By (c),
$\mu_\delta\left(\bigcup_{n\ge 0}Y_n\right) =1$.  Therefore there
exists some $n_0$ \st $\mu_\delta(Y_{n_0})>0$.  Since $Y_n\subset
Y_{n+1}$ for all $n$, we must have $\mu_\delta(Y_n)>0$ for $n \ge
n_0$.

Now let $n$ be \st $\mu_\delta(Y_n\cap U^c)>0$ and let $Y:=Y_n\cap
U^c$. Then $f^j(Y)\cap Y=\es$ for all $j\ge n$.  Now let
$k:=\max\{j:\mu_\delta(f^j(Y)\cap Y)>0\}<n$ and $Y':=f^k(Y)\cap
Y$.  Then $f^{n-k}(Y')\subset U$, so certainly $\mu_\delta(f^j(Y')
\cap Y')=0$ for $j\ge n-k$.  Moreover, by definition of $k$,
$\mu_\delta(f^j(Y') \cap Y')=0$ for $1\le j< n-k$.  Clearly $Y'$
is a wandering set of positive measure, so $\mu_\delta$ is
dissipative.

We now prove the last statement of the lemma.  Suppose that
$\{D_n\}_n$ is the countable basis of the topology of $\C$. Define
$E_n:=\{z\in \C: f^j(z)\notin D_n \hbox{ for all } j\ge 0\}$. If
the statement is not true then $\mu_\delta\left(\bigcup_{n, D_n
\cap \J \neq \emptyset}E_n\right)=1$. So there exists $n$ \st
$\mu_\delta(E_n)>0$.  Hence there is a disk $D$ intersecting $\J$
and a forward invariant set $X$ of positive measure \st $D\cap
X=\es$. Since $f^N(D) \supset \J$ for some $N$, there exists
$Y\subset D$ \st $f^N(Y)= X$. Since $X$ has positive measure, $Y$
must also have positive measure.  Moreover, $f^n(Y)\cap Y=\es$ for
$n\ge N$.

As above, let $k=\max\{j:\mu(f^j(Y)\cap D)>0\}<N$ and $Y':=
f^k(Y)\cap D$. By definition, $\mu_\delta(Y')>0$ and $Y'\subset
D$.  Then $f^{N-k}(Y')\subset X$, so certainly $\mu_\delta(f^j(Y')
\cap Y')=0$ for $j\ge N-k$.  Moreover, by definition of $k$,
$\mu_\delta(f^j(Y') \cap Y')=0$ for $1\le j< N-k$.  Clearly $Y'$
is a wandering set of positive measure, so $\mu_\delta$ is
dissipative.
\end{proof}

When discussing whether $(X,T,m)$ is dissipative or not, it is
important to distinguish between totally dissipative (i.e., there
is no invariant subset $Y \subset X$ of positive measure on which
$(Y,T,m)$ is conservative), or only ``partially'' dissipative,
when such a proper subset exists.

It is easy to see that given a dynamical system $(X, T)$, an
invariant probability measure $\mu$ must be conservative, since, as in the Poincar\'e recurrence Theorem, if
$U\subset X$ is a wandering set of positive measure then
$\mu\left(\bigcup_{n\ge 0}f^{-n}(U)\right)= \sum_{n\ge
0}\mu(f^{-n}(U))=\sum_{n\ge 0}\mu(U)=\infty$: a contradiction.
However, there are systems where $\mu\ll m$ and $\mu$ is an acip,
but $m$ is dissipative (although not totally dissipative). An
example of this is a renormalisable interval map such that the
renormalisation has an acip (supported on a proper subset of the
interval), while Lebesgue measure is dissipative. The following
lemma shows that this does not occur in our setting.

\begin{lemma}
If $\mu_\delta$ is dissipative, then there is no acip, unless the
acip is supported on $\omega(0)$. \label{lem:diss no acip}
\end{lemma}

\begin{proof}
We suppose that $\mu$ is an acip for $\mu_\delta$ such that
$\mu(\omega(0))<1$. We show that $\mu$ is dissipative if and only
if $\mu_\delta$ is dissipative. We use the method of Lemma 8, the
slight difference being that we are concerned with $\mu$, which is
not conformal. But since an invariant probability measure cannot
be dissipative, the lemma follows.

By Lemma 8, we know that $\mu_\delta$ is dissipative if and only
if $f^n(z)\to \omega(0)$ for $\mu_\delta$--a.e. $z\in \J$. This
condition immediately implies that $f^n(z)\to \omega(0)$ for
$\mu$--a.e. $z\in \J$. We now show that $f^n(z)\to \omega(0)$ for
$\mu$--a.e. $z\in \J$ implies that $\mu$ is also dissipative.

Let $U$ be a neighbourhood of $\omega(0)$ with $0<\mu(U)<1$.
Define $Y_n:=\{z\in \C:f^j(z)\in U \hbox{ for all } j\ge n\}$.  By
assumption, $\mu\left(\bigcup_{n\ge 0}Y_n\right) =1$. Therefore
there exists some $n_0$ \st $\mu(Y_{n_0})>0$. Since $Y_n\subset
Y_{n+1}$ for all $n$, we must have $\mu(Y_n)>0$ for $n \ge n_0$.

Now let $n$ be \st $\mu(Y_n\cap U^c)>0$ and let $Y:=Y_n\cap U^c$.
Then $f^j(Y)\cap Y=\es$ for all $j\ge n$.  Let
$k:=\max\{j:\mu(f^j(Y)\cap Y)>0\}<n$.  Let $Y':=f^k(Y)\cap Y$.
Then $f^{n-k}(Y')\subset U$, so certainly $\mu(f^j(Y') \cap Y')=0$
for $j\ge n-k$.  Moreover, by definition of $k$, $\mu(f^j(Y') \cap
Y')=0$ for $1\le j< n-k$.  Clearly $Y'$ is a wandering set of
positive $\mu$--measure, so $\mu$ is dissipative, as required.
\end{proof}

\begin{proof}[Proof of Lemma~\ref{lem:cons}]
Recall that $(\J,f,\mu)$ is dissipative if there is a set $A$ with
$\mu(A)>0$ \st $\mu(f^{-n}(A)\cap A)=0$ for all $n\ge 1$.
Otherwise $(\J,f,\mu)$ is conservative.  Lemma~\ref{lem:diss no
acip} implies that $\mu_\delta$ is conservative.   Hence
Lemma~\ref{lem:diss} implies that $\mu_\delta$--a.e. point has a
dense orbit in $\J$. Therefore, for any open set $U$ intersecting
$\J$ has $\mu_\delta\left(\bigcup_{n\ge 0}f^{-n}(U)\right)=1$.
Suppose that $\mu(\bigcup_{n\ge 0}f^{-n}(W(\theta))) = 0$, so the
complement has positive $\mu$-measure. Then by absolute continuity
$\mu_\delta\left(\left(\bigcup_{n\ge
0}f^{-n}(W(\theta))\right)^c\right)>0$ which is a contradiction to
the previous line. Therefore, $\mu(W(\theta))>0$.
\end{proof}

\section*{Appendix B}

Here we present some lemmas concerning the combinatorial
properties of unimodal maps $f:\R \to \R$ and resulting Koebe
space.

\begin{lemma}\label{lem:basic props2}
The following properties hold:
\begin{itemize}
\item[(a)] $D_{S_k} = f^{S_k}([\z_{k-1}, 0]) = [f^{S_{Q(k)}}(0),
f^{S_k}(0)]$;
\item[(b)] if $J$ is a branch of $f^n$ and $f^n(J)\ni 0$ then $f^n(J) =
D_{S_k}$ for some $k$;
\item[(c)] $D_n = [f^n(0), f^{\beta(n)}(0)]$ for
$\beta(n) = n - \max\{ S_j \ : \ S_j < n\}$ and if $n$ is not a
cutting time, then $D_n \subset D_{\beta(n)}$;
\item[(d)] if $Q(k) \to \infty$, then the lengths $|D_n| \to 0$ as
$n \to \infty$.
\end{itemize}
\end{lemma}

\begin{proof} See \cite{BruinKnea}. \end{proof}

\begin{lemma}\label{lem:no2cpp}
Assume that $Q(k) \to \infty$ satisfies \eqref{eq:strictHofbauer}.
Then there exists $m_0 \in \N$ \st for each non-cutting time $n$,
$D_n$ contains at most one closest precritical point $\z_m$ or
$\hat \z_m$ for $m \ge m_0$.
\end{lemma}

\begin{proof}
From Lemma~\ref{lem:basic props2}(c), we can build, for any
$n$, a nested sequence
\[
D_n \subset D_{\beta(n)} \subset D_{\beta^2(n)} \subset
\cdots\subset D_{\beta^{r-1}(n)} \subset D_{\beta^r(n)},
\]
where $\beta^r(n)$ is a cutting time (say $S_{a(n)-1}$) but
$\beta^{r-1}(n) = S_{a(n)-1} + S_{b(n)-1}$ is not a cutting time.
See Figure~\ref{fig2}. We let $a=a(n),\ b=b(n)$.

By Lemma~\ref{lem:basic props2}(d), $|D_{S_{b-1} + S_{a-1}}| \to
0$ as $b \to \infty$, so there is a neighbourhood $U_a$ of $0$
such that $D_{S_{b-1} + S_{a-1}} \cap U_a = \emptyset$ for all
$b \in \N$ such that $Q(b) > a-1$.
Take $U = \cap_{n,\ a(n) < k_0} U_{a(n)}$ with $k_0$ as
in \eqref{eq:strictHofbauer} and let $m_0$ be minimal such that
$\z_m, \hat \z_m \in U$ for all $m \ge m_0$.

Now take $n$ a non-cutting time, and assume by contradiction that
$D_n$ contains $\z_m$ and $\z_{m+1}$ (or $\hat \z_m$ and $\hat
\z_{m+1}$) for $m \ge m_0$.

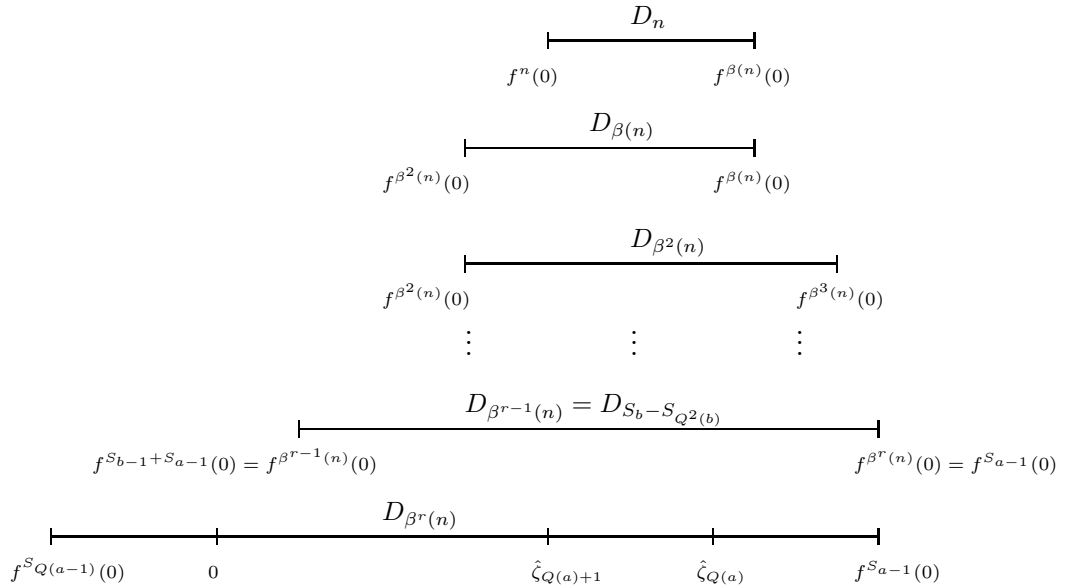
\begin{figure}[ht]
\unitlength=11mm
\begin{picture}(12,7.5)(0,0)
\thinlines \put(1,1){\line(1,0){10}} \put(1,0.9){\line(0,1){0.2}}
\put(0.5,0.5){\tiny $f^{S_{Q(a-1)}}(0)$}
\put(3,0.9){\line(0,1){0.2}} \put(2.9,0.5){\tiny $0$}
\put(11,0.9){\line(0,1){0.2}} \put(10.7,0.5){\tiny
$f^{S_{a-1}}(0)$} \put(9,0.9){\line(0,1){0.2}} \put(8.8,0.5){\tiny
$\hat \z_{Q(a)}$} \put(7,0.9){\line(0,1){0.2}} \put(6.8,0.5){\tiny
$\hat \z_{Q(a)+1}$}

\put(4,2.3){\line(1,0){7}} \put(4,2.2){\line(0,1){0.2}}
\put(1.5,1.8){\tiny $f^{S_{b-1} + S_{a-1}}(0) = f^{\beta^{r-1}(n)}(0)$}
\put(11,2.2){\line(0,1){0.2}} \put(10.7,1.8){\tiny
$f^{\beta^r(n)}(0) = f^{S_{a-1}}(0)$}

\put(6,4.3){\line(1,0){4.5}} \put(6,4.2){\line(0,1){0.2}}
\put(5,3.8){\tiny $f^{\beta^2(n)}(0)$}
\put(10.5,4.2){\line(0,1){0.2}} \put(10,3.8){\tiny
$f^{\beta^3(n)}(0)$}

\put(6,5.7){\line(1,0){3.5}} \put(6,5.6){\line(0,1){0.2}}
\put(5,5.2){\tiny $f^{\beta^2(n)}(0)$}
\put(9.5,5.6){\line(0,1){0.2}} \put(9,5.2){\tiny
$f^{\beta(n)}(0)$}

\put(9.5,7){\line(-1,0){2.5}} \put(9.5,6.9){\line(0,1){0.2}}
\put(9, 6.5){\tiny $f^{\beta(n)}(0)$}
\put(7,6.9){\line(0,1){0.2}} \put(6.5,6.5){\tiny $f^{n}(0)$}

\put(8,7.2){\small $D_n$} \put(7.5,5.9){\small $D_{\beta(n)}$}
 \put(8,4.5){\small $D_{\beta^2(n)}$}
\put(6,2.5){\small $D_{\beta^{r-1}(n)} = D_{S_b - S_{Q^2(b)}}$}
\put(5,1.2){\small $D_{\beta^r(n)}$}

\put(6,3.2){$\vdots$}\put(8,3.2){$\vdots$}\put(10,3.2){$\vdots$}

\end{picture}
\caption{The nested sequence of intervals
$\left(D_{\beta^j(n)}\right)_{j = 0}^r$.} \label{fig2}
\end{figure}

Build the nested sequence $(D_{\beta^j(n)})_{j = 0}^r$ as above
and let $\beta^r(n) = S_{a-1}$ and $\beta^{r-1}(n) = S_{a-1} +
S_{b-1}$. Since $D_n$ intersects $U$, so does $D_{ S_{a-1} +
S_{b-1} }$, and hence $a \ge k_0$.

Since $D_{S_{a-1}+S_{b-1}}$ contains two closest precritical
points and, by Lemma~\ref{lem:basic props}(b), $f^{S_{a-1}}(0) \in
(\z_{Q(a)-1}, \z_{Q(a)}) \cup (\hat \z_{Q(a)}, \hat \z_{Q(a)-1})$,
at least $\z_{Q(a)}$ and $\z_{Q(a)+1}$ (or $\hat \z_{Q(a)}$ and
$\hat \z_{Q(a)+1}$ as in Figure~\ref{fig2})
belong to $D_{S_{a-1} + S_{b-1}}$. It follows
that the first cutting time beyond $S_{a-1} + S_{b-1}$ is $S_b =
S_{b-1}+S_{a-1} + S_{Q(a)}$, so $a = Q(b)$ and
\[
S_{b-1}+S_{a-1}
= S_b - S_{Q(a)} = S_b - S_{Q^2(b)}.
\]
Now apply $f^{S_{Q^2(b)}}$ to the interval $D_{S_b - S_{Q^2(b)}}$;
it maps $f^{S_{b-1} + S_{a-1}}(0) = f^{S_b - S_{Q^2(b)}}(0)$ to $f^{S_b}(0)$,
$\hat \z_{Q^2(b)}$ to $0$ and $\hat \z_{Q^2(b)+1}$ to a point
$\z \in f^{S_{Q^2(b)} - S_{ Q^2(b)+1 } }(0) = f^{-S_{Q(Q^2(b)+1)}}(0)$.
In addition, $\z$ has to be a closest precritical point; otherwise
there is a point $\z' \in (\z,0) \cap f^{-k}(0)$ for some $k \le S_{Q^2(b)}$
that pulls back under $f^{S_{Q^2(b)}}$ to a closest precritical point
strictly between $\hat\z_{Q^2(b)}$ and $\hat\z_{Q^2(b)+1}$.
Thus $\z= \z_{Q(Q^2(b)+1)} \in [f^{S_{Q^2(b)}}(0), 0]$.
The closest precritical point of lowest index in $[f^{S_b}(0), 0]$
is $\z_{Q(b+1)}$ by Lemma~\ref{lem:basic props}(b),
so $Q(b+1) \le Q(Q^2(b)+1)$. But this contradicts assumption
\eqref{eq:strictHofbauer}.
\end{proof}

\begin{lemma}\label{lem:monotone}
Assume that $Q(k) \to \infty$ and that Condition
\eqref{eq:strictHofbauer} holds. Let $m_0$ be as in
Lemma~\ref{lem:no2cpp}. Then if $V \subset \R$ is the largest
neighbourhood of $c_1$ on which $f^{S_k-1}$ is monotone, then
$f^{S_k-1}(V) = [f^{S_{Q^2(k)}}(0), f^{S_{Q(k)}}(0)]$.
\end{lemma}

\begin{proof}
Let $V = [f(w), f(z) ] \owns c_1$ be the largest interval on which
$f^{S_k-1}$ is monotone. By definition of closest precritical
point, $z = \z_{k-1}$, and $f^{S_k-1}(f(\z_{k-1})) =
f^{S_k-S_{k-1}}(0) = f^{S_{Q(k)} }(0)$. The other endpoint $f(w)$
is an image of $w \in i \R$. Decompose $f^{S_k-1} =  f^{S_{Q(k)}}
\circ f^{S_{k-1}-1}$. By Lemma~\ref{lem:basic props}(c), $f^{
S_{k-1}-1 }(c_1) \in (\z_{Q(k)-1}, \z_{Q(k)}) \cup (\hat
\z_{Q(k)}, \hat \z_{Q(k)-1})$; for simplicity assume that
$f^{S_{k-1}-1}(c_1) \in (\z_{Q(k)-1}, \z_{Q(k)})$, as in
Figure~\ref{fig1}. There are two possibilities:

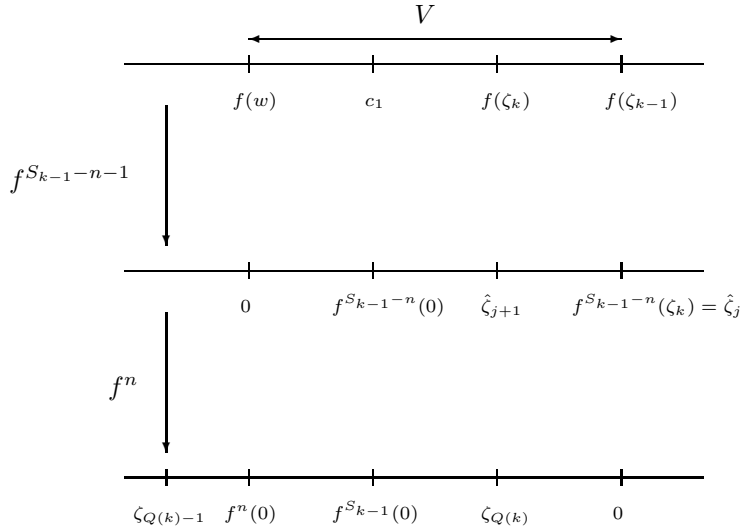
\begin{figure}[ht]
\unitlength=11mm
\begin{picture}(10,7)(0,0)
\thinlines \put(1,6){\line(1,0){7}}
\put(4.5,6.3){\vector(-1,0){2}}\put(4,6.3){\vector(1,0){3}}
\put(4.5,6.5){\small $V$} \put(2.5,5.9){\line(0,1){0.2}}
\put(2.3,5.5){\tiny $f(w)$} \put(4,5.9){\line(0,1){0.2}}
\put(3.9,5.5){\tiny $c_1$} \put(5.5,5.9){\line(0,1){0.2}}
\put(5.3,5.5){\tiny $f(\z_k)$} \put(7,5.9){\line(0,1){0.2}}
\put(6.8,5.5){\tiny $f(\z_{k-1})$}
\put(1.5,5.5){\vector(0,-1){1.7}} \put(-0.4,4.5){\small
$f^{S_{k-1}-n-1}$}

\put(1,3.5){\line(1,0){7}} \put(2.5,3.4){\line(0,1){0.2}}
\put(2.4,3){\tiny $0$} \put(4,3.4){\line(0,1){0.2}}
\put(3.5,3){\tiny $f^{S_{k-1}-n}(0)$}
\put(7,3.4){\line(0,1){0.2}} \put(6.4,3){\tiny
$f^{S_{k-1}-n}(\z_k) = \hat \z_j$} \put(1.5,3){\vector(0,-1){1.7}}
\put(0.8,2){\small $f^n$} \put(5.5,3.4){\line(0,1){0.2}}
\put(5.3,3){\tiny $\hat \z_{j+1}$}

\put(1,1){\line(1,0){7}} \put(1.5,0.9){\line(0,1){0.2}}
\put(1.1,0.5){\tiny $\z_{Q(k)-1}$} \put(2.5,0.9){\line(0,1){0.2}}
\put(2.2,0.5){\tiny $f^n(0)$} \put(4,0.9){\line(0,1){0.2}}
\put(3.5,0.5){\tiny $f^{S_{k-1}}(0)$}
\put(5.5,0.9){\line(0,1){0.2}} \put(5.3,0.5){\tiny $\z_{Q(k)}$}
\put(7,0.9){\line(0,1){0.2}} \put(6.9,0.5){\tiny $0$}

\end{picture}
\caption{The interval $V$ and its images under $f^{S_{k-1}-n-1}$
and $f^{S_{k-1}-1}$.} \label{fig1}
\end{figure}

\begin{itemize}
\item $\z_{Q(k)-1} \in f^{S_{k-1}-1}(V)$. In this case,
$f^{S_{k-1}-1}(f(w)) = \z_{Q(k)-1}$, because otherwise the
interval $f^{S_{k-1}-1}(V)$ is not mapped monotonically for
another $S_{Q(k)}$ iterates, a contradiction.  Therefore,
\[
f^{ S_k-1 }(f(w)) = f^{ S_{Q(k)} }( \z_{Q(k)-1} ) =  f^{S_{Q(k)}-
S_{Q(k)-1}}(0) =  f^{S_{ Q^2(k) }}(0).
\]
\item $\z_{Q(k)-1} \notin f^{S_{k-1}-1}(V)$,
and then there is $n < S_{k-1}$ such  that $f^{S_{k-1}-1}(f(w)) =
f^n(0)$, and $f^{S_{k-1}-n-1}(f(w)) = 0$. It follows that
$f^{S_{k-1}-n}(\z_{k-1})$ is a precritical point. If it is not a
closest precritical point, then pulling back the first closest
precritical point in $(0, f^{S_{k-1}-n}(\z_{k-1}) )$ would give
another closest precritical point between $\z_{k-1}$ and $0$.
Therefore $f^{S_{k-1}-n}(\z_{k-1}) = \z_j$ or $\hat \z_j$ for some
$j < k-1$ (in Figure~\ref{fig1}, this point is $\hat \z_j$), and
$n = S_j$. But pulling back $\z_{Q(k)}$ for $n$ iterates gives
another closest precritical point in $(f^{S_{k-1}-n}(0), \hat
\z_j)$. (This point must be $\hat \z_{j+1}$, which also shows that
$j+1 > Q(k)$, so $j \ge Q(k) \ge m_0$. ) Therefore $D_{S_{k-1}-n}$
contains two closest precritical points $\hat \z_j, \hat \z_{j+1}$
for $j+1 > j \ge m_0$, and as $S_{k-1}-n$ is not a cutting time,
this violates Lemma~\ref{lem:no2cpp}.
\end{itemize}
\end{proof}

\medskip
\noindent
Department of Mathematics\\
University of Surrey\\
Guildford, Surrey, GU2 7XH\\
UK\\
\texttt{h.bruin@surrey.ac.uk}\\
\texttt{http://www.maths.surrey.ac.uk/showstaff?H.Bruin}

\medskip
\noindent
Department of Mathematics\\
University of Surrey\\
Guildford, Surrey, GU2 7XH\\
UK\\
\texttt{m.todd@surrey.ac.uk}\\
\texttt{http://www.maths.surrey.ac.uk/showstaff?M.Todd}

\end{document}